\documentclass[leqno]{article}

\def\diam{\mathop{\rm diam}}
\def\dist{\mathop{\rm dist}}
\def\re{\mathop{\rm Re}}
\def\im{\mathop{\rm Im}}

\begin{document}

\title{Mappings and Spaces, 2}

\author{Stephen W.\ Semmes\thanks{Some of this material was presented
at the November, 2006 meeting of the American Mathematical Society at
the University of Arkansas.  I would like to thank participants of the
meeting and other readers for their comments and suggestions.}}

\date{}

\maketitle

\begin{abstract}
This paper is concerned with analysis on metric spaces in a variety of
settings and with several kinds of structure.
\end{abstract}

\tableofcontents

\section{Basic notions}
\setcounter{equation}{0}

	Let us begin by reviewing some aspects of analysis on metric
spaces.  Let $(M, d(x, y))$ and $(N, \rho(u, v))$ be metric spaces,
let $\alpha$ be a positive real number, and let $C$ be a nonnegative
real number.  A mapping $f : M \to N$ is said to be
\emph{$C$-Lipschitz of order $\alpha$} if
\begin{equation}
        \rho(f(x), f(y)) \le C \, d(x, y)^\alpha
\end{equation}
for every $x, y \in M$.  For example, constant functions are
$0$-Lipschitz of every order and are the only $0$-Lipschitz functions,
and the identity mapping $f(x) = x$ is $1$-Lipschitz of order $1$ as a
mapping from $M$ to $M$.  If $f(x)$ is a continuously-differentiable
real or complex-valued function on the real line, then $f(x)$ is
$C$-Lipschitz of order $1$ if and only if $|f'(x)| \le C$ for every $x
\in {\bf R}$.

         If $f_1$, $f_2$ are complex-valued functions on a metric
space $M$ which are $C_1$, $C_2$-Lipschitz of order $\alpha$ for some
$\alpha > 0$ and $C_1, C_2 \ge 0$, then it is easy to see that $f_1 +
f_2$ is $(C_1 + C_2)$-Lipschitz of order $\alpha$.  If $f$ is a
complex-valued $C$-Lipschitz function of order $\alpha$ on $M$ and $a$
is a complex number, then $a \, f$ is $(|a| \, C)$-Lipschitz function
of order $\alpha$.  If $f_1$, $f_2$ are also bounded, then the product
$f_1 \, f_2$ is Lipschitz of order $\alpha$ too.  The composition of a
Lipschitz mapping of order $\alpha$ and a Lipschitz mapping of order
$\beta$ is Lipschitz of order $\alpha \, \beta$.  Lipschitz mappings
of any order are uniformly continuous.

        For each $p \in M$, $d(x, p) \le d(y, p) + d(x, y)$ for every
$x, y \in M$ by the triangle inequality, and similarly with $x$ and
$y$ interchanged.  This implies that $f_p(x) = d(x, p)$ is a
real-valued $1$-Lipschitz function of order $1$.  If $0 < \alpha \le
1$ and $r$, $t$ are nonnegative real numbers, then $\max(r, t) \le
(r^\alpha + t^\alpha)^{1/\alpha}$, which implies that
\begin{equation}
       r + t \le \max(r, t)^{1-\alpha} \, (r^\alpha + t^\alpha) 
              \le (r^\alpha + t^\alpha)^{1/\alpha},
\end{equation}
or $(r + t)^\alpha \le r^\alpha + t^\alpha$.  It follows that $f_{p,
\alpha} = d(x, p)^\alpha$ is a $1$-Lipschitz function of order
$\alpha$ for every $p \in M$ when $\alpha \le 1$, for the same reasons
as for $\alpha = 1$.

       By contrast, if $f$ is a real or complex-valued Lipschitz
function of order $\alpha > 1$ on the real line, then $f'(x) = 0$ for
every $x \in {\bf R}$, and $f$ is constant.  The same argument works
on Euclidean spaces of any dimension, but there are metric spaces with
nonconstant Lipschitz functions of order $> 1$.  On the Cantor set
there are nontrivial locally constant functions, for instance.  There
are also connected and locally connected snowflake sets with
nonconstant Lipschitz functions of order $> 1$.

        Let $(M, d(x, y))$ be a metric space, and suppose that $a$,
$b$ are real numbers with $a \le b$ and that $p : [a, b] \to M$ is a
continuous curve in $M$.  If $\mathcal{P} = \{t_\ell\}_{\ell = 1}^n$
is a partition of $[a, b]$, which means that
\begin{equation}
        a = t_0 < t_1 < \cdots < t_n = b,
\end{equation}
then put
\begin{equation}
        \lambda(\mathcal{P}) = \sum_{\ell = 1}^n d(p(t_{\ell - 1}), p(t_\ell)).
\end{equation}
We say that $p$ has finite length in $M$ if there is an upper bound
for $\lambda(\mathcal{P})$ over all partitions $\mathcal{P}$ of $[a,
b]$, in which case the length $\lambda$ of $p$ is defined to be the
supremum of $\lambda(\mathcal{P})$.  This is the same as saying that
$p$ has bounded variation when $M$ is ${\bf R}$ or ${\bf C}$.  If $p :
[a, b] \to M$ is a continuous curve of length $\lambda$, then $d(p(a),
p(b)) \le \lambda$.  The restriction of $p$ to any subinterval of $[a,
b]$ is a continuous curve with length $\le \lambda$, and hence the
diameter of $p([a, b])$ is $\le \lambda$.  Note that a continuous
curve has length equal to $0$ if and only if it is constant.

        Suppose that $\mathcal{P}$, $\mathcal{P}'$ are partitions of
$[a, b]$ and that $\mathcal{P}'$ is a refinement of $\mathcal{P}$,
which is to say that each term in $\mathcal{P}$ is also in
$\mathcal{P}'$.  Using the triangle inequality, one can check that
$\lambda(\mathcal{P}) \le \lambda(\mathcal{P}')$.  As a consequence,
it suffices to use partitions of $[a, b]$ that contain a fixed element
$r \in [a, b]$ in order to determine the length of $p$.  This implies
that the length of $p$ on $[a, b]$ is equal to the sum of the lengths
of $p$ on $[a, r]$ and on $[r, b]$ for each $r \in [a, b]$.  If $p :
[a, b] \to M$ is $C$-Lipschitz of order $1$, then $p$ has finite
length $\le C \, (b - a)$.  Conversely, a continuous path of finite
length $\lambda$ can be reparameterized to get a $1$-Lipschitz curve
on an interval of length equal to $\lambda$.  This basically uses the
arc-length parameterization of $p$.

        If $p : [a, b] \to M$ is a continuous path of length $\lambda$
and $f$ is a $C$-Lipschitz complex-valued function of order $1$ on
$M$, then $f \circ p$ is a function of bounded variation on $[a, b]$
of total variation $\le C \, \lambda$.  If $f$ is locally
$C$-Lipschitz of order $1$, then $f \circ p$ is still a function of
bounded variation on $[a, b]$ of total variation $\le C \, \lambda$,
since one can use partitions of $[a, b]$ with small mesh size by
passing to suitable refinements.  If $f$ is locally
$\epsilon$-Lipschitz of order $1$ for each $\epsilon > 0$, then $f
\circ p$ has total variation equal to $0$, and $f \circ p$ is constant
on $[a, b]$.  If $f$ is Lipschitz of order $> 1$, then $f$ is locally
$\epsilon$-Lipschitz of order $1$ for each $\epsilon > 0$.

         If there is a $k \ge 1$ such that every $x, y \in M$ can be
connected by a continuous path of length $\le k \, d(x, y)$, and if
$f$ is a locally $C$-Lipschitz complex-valued function of order $1$ on
$M$, then $f$ is $(k \, C)$-Lipschitz of order $1$.  This property
holds with $k = 1$ when $M$ is a convex set in Euclidean space or a
normed vector space more generally, since every pair of elements of
$M$ can be connected by a line segment of length equal to the distance
between $x$ and $y$.  This property also holds for some fractals like
the Sierpinski gasket and carpet, for suitable $k > 1$.  A connected
open set $U$ in a normed vector space satisfies this property locally
with $k = 1$, and every $x, y \in U$ can be connected by a curve of
finite length, but the relationship between the lengths of the paths
and the distances between $x$ and $y$ may be complicated.  Similarly,
a connected embedded smooth submanifold of ${\bf R}^n$ has this
property locally with $k$ arbitrarily close to $1$, but otherwise the
ambient Euclidean distance may be much smaller than the intrinsic
distance on the submanifold, depending on the situation.

\section{Generalized pseudomanifold spaces}
\setcounter{equation}{0}

	As on p148 of \cite{spa}, an $n$-dimensional
\emph{pseudomanifold} $M$ is a simplicial complex such that every
simplex in $M$ is contained in an $n$-dimensional simplex in $M$, and
every $(n - 1)$-dimensional simplex in $M$ is a face of one or two
$n$-dimensional simplices in $M$.  It is customary to ask also that
$M$ satisfy the connectedness condition that for every pair of
$n$-dimensional simplices $\sigma$, $\sigma'$ in $M$ there is a
sequence $\sigma_1, \ldots, \sigma_l$ of $n$-dimensional simplices in
$M$ such that $\sigma_1 = \sigma$, $\sigma_l = \sigma'$, and
$\sigma_i$, $\sigma_{i + 1}$ are adjacent when $1 \le i < l$ in the
sense that $\sigma_i \cap \sigma_{i + 1}$ is an $(n - 1)$-dimensional
simplex in $M$.  In particular, this implies that $M$ is connected as
a topological space, since simplices are connected sets.  The
\emph{boundary} $\partial M$ of $M$ consists of the $(n -
1)$-dimensional simplices in $M$ which are faces of exactly one
$n$-dimensional simplex in $M$.  The aforementioned connectedness
condition implies that the \emph{interior} $M \backslash \partial M$
of $M$ is connected too.

	Suppose that $p \in M$ is in the interior of an
$n$-dimensional simplex $\sigma$ in $M$, or in the interior of an $(n
- 1)$-dimensional simplex $\tau$ in $M$ contained in two distinct
$n$-dimensional simplices $\sigma'$, $\sigma''$ in $M$.  In the first
case there is a neighborhood of $p$ in $M$ contained in $\sigma$, and
in the second case there is a neighborhood of $p$ in $M$ contained in
$\sigma' \cup \sigma''$.  In both cases there is a neighborhood of $p$
in $M$ which is homeomorphic to the open unit ball in ${\bf R}^n$.  If
$p \in M$ lies in the interior of an $(n - 1)$-dimensional simplex
$\tau$ in the boundary of $M$, which means that $\tau$ is contained in
exactly one $n$-dimensional simplex $\sigma$ in $M$, then there is a
neighborhood of $p$ in $M$ contained in $\sigma$.  In this event $M$
looks exactly like an $n$-dimensional manifold with boundary at $p$.

	If $p \in M$ is contained in a $k$-dimensional simplex with $k
\le n - 2$, then the local behavior of $M$ at $p$ may be more
complicated.  The set $E$ of these potentially singular points in $M$
is relatively small in $M$, and in particular the connectedness
condition for $M$ implies that $M \backslash E$ is connected.

	An \emph{orientation} on $M$ consists of an orientation on
every $n$-dimensional simplex in $M$, with the compatibility condition
that for every $(n - 1)$-dimensional simplex $\tau$ in $M$ which is
contained in two adjacent $n$-dimensional simplices $\sigma'$,
$\sigma''$, the orientations on $\tau$ induced by those on $\sigma'$,
$\sigma''$ are the same.  If there is an orientation on $M$, then $M$
is said to be \emph{orientable}, as usual.  Because of the
compatibility condition, a choice of orientation on an $n$-dimensional
simplex in $M$ determines orientations on the adjacent $n$-dimensional
simplices, and hence on all of $M$ by the connectedness condition when
$M$ is orientable.  One can try to get an orientation on $M$ by
starting with an orientation on an $n$-dimensional simplex and using
induced orientations on adjacent simplices, etc.  This works when
every $n$-dimensional simplex only gets one orientation in this way.

	It seems natural to consider more general situations with
analogous features.  Let $n$ be a positive integer, let $(M, d(x, y))$
be a separable metric space, and let $A$, $B$ be closed subsets of $M$
such that $M \backslash (A \cup B)$ is dense in $M$.  Roughly
speaking, the singularities of $M$ ought to be contained in $A$, and
the boundary of $M$ ought to be contained in $B$.  Suppose in
particular that the topological dimension of $M$ is equal to $n$, that
the topological dimension of $A$ is less than or equal to $n - 2$, and
that the topological dimension of $B$ is less than or equal to $n -
1$.  Let us ask that $M \backslash (A \cup B)$ be connected as well,
which implies that $M$ is connected.

	As a basic scenario, suppose that $M \backslash A$ is an
$n$-dimensional $C^\infty$ manifold with boundary $B \backslash A$,
where the manifold structure is compatible with the topology
determined by the metric on $M$.  Suppose also that $M \backslash A$
is equipped with a $C^\infty$ Riemannian metric, which leads to a
Riemannian distance function on $M \backslash A$ by minimizing the
lengths of paths in the usual way.  A standard local compatibility
condition asks that there be a $C > 0$ and, for every $p \in M
\backslash A$, an open set $U(p) \subseteq M \backslash A$ such that
$p \in U(p)$ and the Riemannian distance and $d(x, y)$ are each
bounded by $C$ times the other on $U(p)$.  This implies that the
lengths of curves in $M \backslash A$ associated to the Riemannian
structure are comparable to those defined using $d(x, y)$, and that
Riemannian volumes of subsets of $M \backslash A$ are comparable to
$n$-dimensional Hausdorff measure defined using $d(x, y)$.
Comparability of distances associated to the Riemannian structure and
$d(x, y)$ globally on $M \backslash A$ would be an interesting
additional condition.

	In any case, Hausdorff measures and dimensions can be defined
for subsets of $M$ using $d(x, y)$, and the topological dimension of a
set is automatically less than or equal to its Hausdorff dimension.
More precisely, if a set has $l$-dimensional Hausdorff measure equal
to $0$, then its topological dimension is strictly less than $l$.  To
say that the singular set $A$ is relatively small in $M$, one can
consider stronger conditions in terms of Hausdorff measure.  One might
ask that the $(n - 1)$-dimensional Hausdorff measure of $A$ be equal
to $0$, or that that the Hausdorff dimension of $A$ be strictly less
than $n - 1$, or that the $(n - 2)$-dimensional Hausdorff measure of $A$
be finite, perhaps at least locally.  Depending on the circumstances,
one might consider even more restrictive conditions on $A$, associated
to smaller dimensions.

	One might also consider more complicated types of boundaries.
Instead of smoothness, one might consider accessibility conditions,
for instance.  Concerning the size of $B$, one might ask that the
$n$-dimensional Hausdorff measure of $B$ be equal to $0$, or that the
Hausdorff dimension of $B$ be strictly less than $n$, or that the $(n
- 1)$-dimensional Hausdorff measure of $B$ be finite, at least
locally.  It may be that $M$ is equipped with a nonnegative Borel
measure $\mu$, and one might ask that $\mu(A) = \mu(B) = 0$.  If $B
\ne \emptyset$, then one might ask that $B$ be at least $(n -
1)$-dimensional in various ways.

	Let $B(x, r) = \{y \in M : d(x, y) < r \}$ be the open ball in
$M$ with center $x \in M$ and radius $r > 0$.  As in \cite{c-w-1,
c-w-2}, $M$ is \emph{doubling} if there is a $C > 0$ such that for
every $x \in M$ and $r > 0$ there are $x_1, \ldots, x_l \in M$ with $l
\le C$ and
\begin{equation}
	B(x, 2 r) \subseteq \bigcup_{i = 1}^l B(x_i, r).
\end{equation}
Similarly, a nonnegative Borel measure $\mu$ on $M$ is a
\emph{doubling measure} if the $\mu$-measure of every ball in $M$ is
finite, and if there is a $C' > 0$ such that
\begin{equation}
	\mu(B(x, 2 r)) \le C' \, \mu(B(x, r))
\end{equation}
for every $x \in M$ and $r > 0$.  One can show that $M$ is doubling
when there is a nonzero doubling measure on $M$.

	A set $E \subseteq M$ is said to be \emph{porous} if there is
a $c > 0$ such that for every $x \in M$ and $r > 0$, with $r$ less
than or equal to the diameter of $M$ when $M$ is bounded, there is a
$z \in M$ for which $d(x, z) < r$ and $B(z, c \, r) \cap E =
\emptyset$.  The closure of a porous set is clearly porous, with the
same constant $c$, and one can check that the union of two porous sets
is porous.  Equivalently, $E \subseteq M$ is porous if there is a $c'
> 0$ such that for every $x \in E$ and $r > 0$ with $r \le \diam M$
when $M$ is bounded, there is a $z \in M$ which satisfies $d(x, z) <
r$ and $B(z, c' \, r) \cap E = \emptyset$.  Consequently, if $E
\subseteq Y \subseteq M$ and $E$ is porous as a set in $Y$, then $E$
is porous as a set in $M$.

	For a generalized pseudomanifold space $M$ with singular set
$A$ and boundary $B$, a doubling condition for $M$ and perhaps a
nonnegative measure $\mu$ on $M$ as well as porosity conditions on
$A$, $B$ in $M$ can be quite appropriate.  One of the nice features of
doubling measures is that there is a version of the Lebesgue density
theorem, so that sets of positive measure have points of density.  It
follows that porous sets have measure $0$ with respect to doubling
measures, because they cannot have points of density in this case.  In
${\bf R}^n$, one can show that porous sets have Hausdorff dimension
strictly less than $n$, and there are related results for other spaces
depending on the circumstances.

	An appealing quantitative local connectedness property for the
regular part of pseudomanifold space would be the same as for uniform
domains \cite{mar-s}, although this would have to be relaxed for some
singularities.  Higher-order versions \cite{als, a-v-1, a-v-2, h-y} can
be quite interesting too, especially for more precise information
about the structure of the singularities.  These are refinements of
local connectedness conditions commonly studied in geometric topology.
See \cite{b-h-k, h2, mar} for more information about uniform domains.

	One can also consider mixtures of local regularity and
relatively small singular sets for functions on pseudomanifold spaces
and mappings between them.  A standard argument would combine
estimates for the size of a singular set $A$ and continuity properties
of a mapping $f$ on $A$ to estimate the size of $f(A)$.  If $f(A)$ is
sufficiently small, then one may be able to avoid the singularities
and work on the regular part, as in the study of degrees of mappings
and other topological properties.  In particular, one could use the
implicit function theorem on the regular part.  For that matter, a
function might have its own singularities, and the pseudomanifold
point of view can be a convenient way to adapt the geometry of a space
to the behavior of a function on it.

	The singular set in \cite{l2} is a nontrivial connected set in
a $2$-dimensional space.  This suggests some variants of some of the
questions in \cite{h-s}, e.g., concerning the existence of bilipschitz
coordinates for $2$-dimensional generalized pseudomanifold spaces,
under suitable conditions.  Specifically, one can consider situations
in which the singular set is uniformly disconnected.  This would be a
version of the hypothesis that the singular set have topological
codimension at least $2$.  However, it would allow the Hausdorff
codimension of the singular set to be arbitrarily small.

	One can consider similar notions for other types of objects,
like weights.  Let $(M, d(x, y))$ be a metric space, and let $w(x)$ be
a nonnegative extended real-valued function on $M$, i.e., $0 \le w(x)
\le + \infty$ for $x \in M$.  Let $U$ be the set of $x \in M$ for
which there are positive real numbers $\epsilon$, $k$, $r$ such that
$\epsilon \le w \le k$ on $B(x, r)$, and note that $U$ is
automatically an open set in $M$.  As a basic class of weights on $M$,
one can consider the $w$'s for which $U$ is dense in $M$ too.  One
might also ask that $w$ be continuous on $U$, or on all of $M$ using
the standard topology for the extended real numbers.

	It is easy to formulate quantitative scale-invariant
regularity conditions for weights on $M$.  For instance, suppose that
there are positive real numbers $c_1$, $c_2$ such that for every $x
\in M$ and $r > 0$ with $r \le \diam M$ when $M$ is bounded there is a
$z \in M$ which satisfies $d(x, z) < r$, $0 < w < +\infty$ on $B(z,
c_1 \, r)$, and $w(y) \le c_2 \, w(y')$ for $y, y' \in B(z, c_1 \,
r)$.  In particular, this implies that $M \backslash U$ is porous in
$M$ with constant $c_1$.  Alternatively, one might start with a
nonempty porous set $A \subseteq M$, and ask that $0 < w < +\infty$ on
$M \backslash A$, and that for every $l \ge 1$ there be a $C(l) \ge 1$
such that $w(y) \le C(l) \, w(y')$ when $y, y' \in M \backslash A$ and
\begin{equation}
\label{d(y, y') le l min(dist(y, A), dist(y', A))}
	d(y, y') \le l \, \min(\dist(y, A), \dist(y', A)).
\end{equation}
Here $\dist(x, A)$ is the infimum of $d(x, u)$ over $u \in A$, as
usual, and $w(x) = \dist(x, A)^\alpha$ has this property for every
$\alpha \in {\bf R}$, because (\ref{d(y, y') le l min(dist(y, A),
dist(y', A))}) implies that the distances from $y$, $y'$ to $A$ are
each less than or equal to $l + 1$ times the other.

	Now suppose that $f$ is an extended real-valued function on
$M$, and let $V$ be the open set of $x \in M$ for which there is an $r
> 0$ such that $f(B(x, r))$ is a bounded set in ${\bf R}$.  We can
begin by asking that $V$ be dense in $M$, and perhaps that $f$ be
continuous on $V$ or on all of $M$.  As a quantitative scale-invariant
condition, we can ask that there be $c, C > 0$ such that for every $x
\in M$ and $r > 0$ with $r \le \diam M$ when $M$ is bounded there is a
$z \in M$ with $d(x, z) < r$, $f(z) \in {\bf R}$, and $|f(y) - f(z)|
\le C$ when $d(y, z) < c \, r$, which implies that $M \backslash V$ is
porous with constant $c$.  As a stronger condition, we can ask that
there be a nonempty porous set $A \subseteq M$ such that $f$ is
real-valued on $M \backslash A$ and $|f(y) - f(y')|$ is bounded when
$y, y' \in M \backslash A$ satisfy (\ref{d(y, y') le l min(dist(y, A),
dist(y', A))}), with a bound that depends on $l$.  These are variants
of ``bounded mean oscillation'' which correspond to logarithms of
weights as in the previous paragraphs.

	A \emph{quasisymmetric mapping} \cite{tuk-v} from one metric
space to another approximately preserves relative distances.  If
$M_1$, $M_2$ are metric spaces, $E \subseteq M_1$ is a porous set, and
$\phi : M_1 \to M_2$ is quasisymmetric, then one can show that
$\phi(E)$ is porous in $M_2$.  As in \cite{v2}, if $M_1 = M_2 = {\bf
R}^n$ with the standard metric, $E \subseteq {\bf R}^n$ is porous, and
$\phi : E \to {\bf R}^n$ is a quasisymmetric mapping, then $\phi(E)$
is porous in ${\bf R}^n$.  One can reformulate (\ref{d(y, y') le l
min(dist(y, A), dist(y', A))}) as saying that $d(y, y') \le l
\min(d(y, a), d(y', a))$ for every $a \in A$, which is preserved by a
quasisymmetric mapping, if we are allowed to replace $l$ by a positive
real number depending on $l$ and the quasisymmetry condition for the
mapping.  The corresponding local regularity conditions for functions
and weights are therefore preserved by quasisymmetric mappings too.

\section{Complex-analytic metric spaces}
\setcounter{equation}{0}

	What might one mean by a ``complex-analytic metric space''?
Certainly ${\bf C}^n$ with the standard Euclidean metric ought to be
an example, as well as domains in ${\bf C}^n$ and smooth complex
manifolds equipped with suitable geometries, etc.

	For a nonstandard example, fix an integer $n \ge 2$, and let
$\Sigma_n$ be the unit sphere in ${\bf C}^n$.  Thus
\begin{equation}
	\Sigma_n = \{z \in {\bf C}^n : |z| = 1\},
\end{equation}
where $|z| = \Big(\sum_{j = 1}^n |z_j|^2 \Big)^{1/2}$ for $z = (z_1,
\ldots, z_n) \in {\bf C}^n$, as usual.  We can think of $\Sigma_n$ as
a real smooth hypersurface in ${\bf C}^n$, whose tangent space at $z
\in \Sigma_n$ is
\begin{equation}
	T_z \, \Sigma_n = \Big\{v \in {\bf C}^n : 
		\re \sum_{j = 1}^n v_j \, \overline{z_j} = 0 \Big\}.
\end{equation}
Here $\re a$ denotes the real part of a complex number $a$, and
$\overline{a}$ is its complex conjugate.  Put
\begin{equation}
	CT_z \, \Sigma_n = \Big\{v \in {\bf C}^n :
			\sum_{j = 1}^n v_j \, \overline{z_j} = 0 \Big\},
\end{equation}
which is a complex-linear subspace of ${\bf C}^n$ contained in $T_z \,
\Sigma_n$.  If $\im a$ is the imaginary part of a complex number $a$,
then $CT_z \, \Sigma_n$ consists of the $v \in T_z \, \Sigma_n$ such
that $\im \sum_{j = 1}^n v_j \, \overline{z_j} = 0$.  This shows that
$CT_z \, \Sigma_n$ has real codimension $1$ in $T_z \, \Sigma_n$,
which has real codimension $1$ in ${\bf C}^n$.  It is well known that
every pair of elements of $\Sigma_n$ can be connected by a smooth path
$p(t)$ in $\Sigma_n$ whose derivative $\dot p(t)$ is contained in
$CT_{p(t)} \, \Sigma_n$ for every $t$ in the interval on which $p(t)$
is defined.  A metric on $\Sigma_n$ can be defined using the infimum
of the lengths of these paths with a fixed pair of endpoints in
$\Sigma_n$.  With respect to this sub-Riemannian geometry on
$\Sigma_n$, $CT_z \, \Sigma_n$ is the appropriate tangent space for
$\Sigma_n$ at $z \in \Sigma_n$.  By construction, $CT_z \, \Sigma_n$
is also a complex vector space in a natural way.

	This sub-Riemannian geometry on $\Sigma_n$ is compatible with
the usual topology, but the corresponding Hausdorff dimension is $2 \,
n$.

	On a complex manifold $M$, there is a decomposition of
exterior differentiation $d$ into the sum of $\partial$ and
$\overline{\partial}$.  By definition, a complex-valued function $f$
on an open set $U \subseteq M$ is holomorphic if $\overline{\partial}
f = 0$ on $U$.  On $\Sigma_n$, there is an analogous operator
$\overline{\partial}_b$ based on the complex subspaces of the tangent
spaces.  The $\overline{\partial}_b$ operator is the appropriate
$\overline{\partial}$ operator on $\Sigma_n$ with respect to the
sub-Riemannian geometry.

	Similar remarks can be applied to other Cauchy--Riemann
manifolds with compatible sub-Riemannian geometries.  In order to get
a complex-analytic metric space, one ought to have complex structures
on the subspaces of the tangent spaces that determine the
sub-Riemannian structure.  Otherwise, one might have ``Cauchy--Riemann
sub-Riemannian spaces'', with complex structures on subspaces of the
subspaces of the tangent spaces that determine the sub-Riemannian
structure.

	In general, one might ask that a complex-analytic metric space
have some sort of tangent spaces, perhaps almost everywhere, and
complex structures on these tangent spaces.  Some nontrivial
holomorphic functions would be nice too.

	Of course, there has been a lot of work over the years
concerning abstract versions of holomorphic functions on complex
spaces, often in terms of algebras of continuous functions on
topological spaces.  An advantage of metric spaces is that there are
special classes of functions, like Lipschitz functions and Sobolev
spaces when the metric space is equipped with a metric, which are
relevant for differentiation and other aspects of analysis.  The
definition of the tangent spaces of the metric space would normally
involve some sort of regular functions and their derivatives.

        Geometric measure theory deals extensively with
differentiation and tangent spaces for sets that may not be smooth.
See \cite{alr, alm, har, har-k, har-l, har-s, k, ll1, ll2, sh, si1,
si2}, for instance, concerning holomorphic chains as currents.

	For a metric space equipped with a doubling measure and for
which there are suitable versions of Poincar\'e inequalities, Cheeger
\cite{che} has shown that there are versions of classical results on
differentiability almost everywhere.  This is a very interesting
setting in which to consider $\overline{\partial}$ operators.

	In particular, one might do this for a space $X$ which is a
Cartesian product of an even number of spaces like those described by
Laakso \cite{l1} and intervals.  If $L$ is a Laakso space, then there
is a natural projection from the product of a Cantor set $C$ and the
unit interval $I$ onto $L$, and another projection from $L$ onto $I$.
The composition of these two mappings is the usual coordinate
projection from $C \times I$ onto $I$.  If a complex structure is
defined on $X$ in a compatible way, then one can use these projections
onto intervals to get nontrivial holomorphic functions on $X$.  One
can jazz this up a bit using branching.

	One can also look at holomorphic mappings between
complex-analytic metric spaces, e.g., nontrivial analytic disks.  It
is well known that any holomorphic mapping from a disc into the unit
sphere $\Sigma_n$ in ${\bf C}^n$ is constant.  There are plenty of
analytic disks in a product of Laakso spaces and intervals when the
complex structure on the product satisfies suitable compatibility
conditions.

	For that matter, one could view the product of a Cantor set
and ${\bf C}^n$ as a complex-analytic metric space, in which only the
complex structure in the ${\bf C}^n$ directions is employed in the
product.  Thus a holomorphic function on the product would be
holomorphic on each copy of ${\bf C}^n$, and an analytic disk in the
product would be an analytic disk in one of the copies of ${\bf C}^n$.
More precisely, the Cantor set would be treated as not contributing to
the tangent space of the product or the complex structure.  This is
consistent with the failure of differentiability theorems for
Lipschitz functions on Cantor sets, although one might say instead
that the derivative is equal to $0$.

	Alternatively, let $E$ be a closed set in ${\bf C}^n$, and
suppose that $f : E \to {\bf C}$ is continuously differentiable in the
sense of Whitney.  This means that for each $p \in E$ there is a
real-linear mapping $df_p : {\bf C}^n \to {\bf C}$ which is continuous
as a function of $p$ and satisfies
\begin{equation}
	f(z) = f(p) + df_p(z - p) + o(1)
\end{equation}
uniformly on compact subsets of $E$.  The restriction to $E$ of a
continuously-differentiable function on ${\bf C}^n$ automatically has
this feature, using the ordinary differential of $f$ at $p \in E$.
Conversely, a function $f$ on $E$ with this property has an extension
to a continuously differentiable function on ${\bf C}^n$ whose
differential at $p \in E$ is equal to $df_p$, by Whitney's extension
theorem.  If $f$ is a continuously-differentiable function on $E$ in
the sense of Whitney, then $\overline{\partial} f_p$ can be defined
using $df_p$ in the usual way, and $\overline{\partial} f_p = 0$ for
every $p \in E$ when $f$ is the restriction to $E$ of a holomorphic
function on an open set $U \subseteq {\bf C}^n$ containing $E$.

	If $E$ is a Cantor set or a snowflake, then there are
nontrivial functions $f$ on $E$ which are continuously-differentiable
in Whitney's sense with $df_p = 0$ for every $p \in E$.  At the
opposite extreme, suppose that $E = {\bf R}^n \subseteq {\bf C}^n$ and
$f : {\bf R}^n \to {\bf C}$ is continuously-differentiable as a
function on ${\bf R}^n$.  The differential of $f$ at $p \in {\bf R}^n$
is therefore defined as a real-linear mapping from ${\bf R}^n$ to
${\bf C}$, which has a unique extension to a complex-linear mapping
from ${\bf C}^n$ to ${\bf C}$.  If we use this extension as $df_p$,
then $\overline{\partial} f_p = 0$.

	A basic issue about complex-analytic metric spaces is the
strength of the $\overline{\partial}$ operator, starting with the
question of whether $|\overline{\partial} f|$ is roughly like $|d f|$
when $f$ is real-valued.  This is an elementary feature of the
classical case, and there is an analogous statement for
Cauchy--Riemann spaces in terms of the tangential part of the
differential.  However, this does not say much about the strength of
the $\overline{\partial}$ operator applied to complex-valued
functions, since there are standard local regularity results for
holomorphic functions on ${\bf C}^n$ while the boundary values of
holomorphic functions on the unit ball automatically satisfy the
tangential Cauchy--Riemann equations on the unit sphere but do not
have to be smooth.

	If $f$ is a nice complex-valued function with compact support
on ${\bf C}^n$ and $1 < p < \infty$, then
\begin{equation}
	\int_{{\bf C}^n} |\partial f(z)|^p \, dz
	 \le A(p, n) \, \int_{{\bf C}^n} |\overline{\partial} f(z)|^p \, dz,
\end{equation}
where $A(p, n) > 0$ depends only on $p$ and $n$.  This follows from
well-known results in harmonic analysis, and there are similar
estimates for other norms and spaces of functions.  These matters have
also been studied extensively for domains in ${\bf C}^n$, their
boundaries, and other complex manifolds and Cauchy--Riemann spaces,
with additional terms or boundary conditions, etc., according to the
situation.  Properties like these are of interest for complex-analytic
metric spaces in general, as well as the relationship with a suitable
Laplace operator and subharmonicity.

	The classical theory of quasiconformal mappings in the plane
deals exactly with the Beltrami operators $\overline{\partial}_\mu =
\overline{\partial} - \mu \, \partial$ associated to a perturbation of
the standard complex structure.  The quasiconformality condition
$\|\mu\|_\infty < 1$ ensures that $|\partial_\mu f|$ is comparable to
$|df|$ when $f$ is real-valued.  Moreover, it leads to $L^2$ estimates
for the gradient, and $L^p$ estimates when $p$ is sufficiently close
to $2$.  If a function is holomorphic with respect to
$\overline{\partial}_\mu$, then it can be expressed as the composition
of an ordinary holomorphic function with a quasiconformal mapping with
dilatation $\mu$.

	It can be easier to make sense of the size $|df|$ of the
differential of a function $f$ on a metric space than the differential
$df$, and it may be easier in some situations to make sense of
something like $|\overline{\partial} f|$ than $\overline{\partial} f$.
There could also be a decomposition of $|df|^2$ into a sum of parts
corresponding to $|\partial f|^2$ and $|\overline{\partial} f|^2$,
analogous to the usual decomposition of $d$ into the sum of $\partial$
and $\overline{\partial}$.  One might look at this on the Sierpinski
gasket in connection with ``analysis on fractals'' in the sense of
\cite{kig, sr2, sr3}, for instance.  The underlying local model for
this is the fact that a real-affine function on the plane is uniquely
determined by its values on the vertices of a triangle, and the
decomposition of a real-linear function into parts that are
complex-linear and conjugate-linear.  By contrast, this may not work
as well for squares and Sierpinski carpets.

	Let $(M, d(x, y))$ and $(N, \rho(u, v))$ be metric spaces.  A
mapping $f : M \to N$ is said to be Lipschitz if it is Lipschitz of
order $1$, and it is a bilipschitz embedding of $M$ into $N$ if
$\rho(f(x), f(y))$ is bounded from above and below by constant
multiples of $d(x, y)$ for every $x, y \in M$.  For example, the
standard embedding of the unit sphere $\Sigma_n$ into ${\bf C}^n$ is
bilipschitz with respect to the ordinary Euclidean metric on ${\bf
C}^n$ and the induced Riemannian metric on $\Sigma_n$.  However, this
mapping is Lipschitz and not bilipschitz when one uses the
sub-Riemannian geometry on $\Sigma_n$ associated to the complex
subspaces of the tangent spaces.  There are probably a lot of
subtleties involved with embeddings of complex-analytic metric spaces.

	Note that the boundary values of a holomorphic function on the
unit ball in ${\bf C}^2$ could be considered as a quasiregular mapping
from $\Sigma_2$ with the usual sub-Riemmanian structure into the
complex numbers.  Similalry, the standard projection from the product
of a Laakso space and an interval or another Laakso space to the
complex numbers could also be considered quasiregular.  In these
examples, the tangent spaces of the domain and range have the same
dimension, and quasiregularity can be formulated in terms of the
differentials of the mappings as linear transformations between the
corresponding tangent spaces.  In the first example, the Hausdorff
dimension of the domain is strictly larger than the topological
dimension, which is strictly larger than the dimension of the tangent
spaces.  In the second example, the Hausdorff dimension of the domain
is strictly larger than the topological dimension, which is equal to
the dimension of the tangent spaces.  Even for variants of Laakso's
construction using Cantor sets with Hausdorff dimension $0$ so that
the Hausdorff and topological dimensions of the resulting spaces would
be the same, the Hausdorff measure would not be $\sigma$-finite in the
topological dimension.  The fibers of the mapping are at least totally
disconnected in the second example, if not discrete.  Compare with
\cite{h-hol}.

\section{Clifford holomorphic functions}
\setcounter{equation}{0}

	Let $n$ be a positive integer, and let $\mathcal{C}(n)$ be the
Clifford algebra over the real numbers ${\bf R}$ with $n$ generators
$e_1, \ldots, e_n$.  By definition, $\mathcal{C}(n)$ is an associative
algebra with a nonzero multiplicative identity element.  Thus
$\mathcal{C}(n)$ contains a copy of ${\bf R}$, and the real number $1$
can be identified with the multiplicative identity element of
$\mathcal{C}(n)$.  The generators $e_1, \ldots, e_n$ of
$\mathcal{C}(n)$ satisfy the relations $e_l^2 = -1$ for $l = 1,
\ldots, n$, and $e_q \, e_p = - e_p \, e_q$ when $1 \le p, q \le n$
and $p \ne q$.  For $I = \{l_1, \ldots, l_r\}$, $1 \le l_1 < l_2 <
\cdots < l_r \le n$, let $e_I$ be the element of $\mathcal{C}(n)$
defined by
\begin{equation}
	e_I = e_{l_1} e_{l_2} \cdots e_{l_r}.
\end{equation}
We can include $I = \emptyset$ by putting $e_\emptyset = 1$.  The
$2^n$ elements $e_I$ of $\mathcal{C}(n)$, where $I$ runs through all
subsets of $\{1, \ldots, n\}$, forms a basis for $\mathcal{C}(n)$ as a
vector space over ${\bf R}$.

	Actually, one can think of $\mathcal{C}(n)$ as being equal to
${\bf R}$ when $n = 0$.  When $n = 1$, $\mathcal{C}(n)$ is equivalent
to the complex numbers ${\bf C}$, with the one generator $e_1$
corresponding to the complex number $i$.  When $n = 2$,
$\mathcal{C}(n)$ is equivalent to the quaternions ${\bf H}$.  Normally
one might represent $x \in {\bf H}$ as
\begin{equation}
	x = x_1 + x_2 \, i + x_3 \, j + x_4 \, k,
\end{equation}
where $i^2 = j^2 = k^2 = -1$, $k = i \, j$, and $j \, i = - k$, which
yield $i \, k = -j = - k \, i$ and $j \, k = i = - k \, j$.  For the
identification with $\mathcal{C}(2)$, $i, j \in {\bf H}$ correspond to
the two generators $e_1, e_2 \in \mathcal{C}(2)$, and $k \in {\bf H}$
corresponds to their product $e_1 \, e_2$.

	If $v = (v_1, \ldots, v_n) \in {\bf R}^n$, then we can
associate to $v$ the element
\begin{equation}
	\widehat{v} = v_1 \, e_1 + \cdots + v_n \, e_n
\end{equation}
of $\mathcal{C}(n)$.  This defines a linear embedding of ${\bf R}^n$
into $\mathcal{C}(n)$ such that
\begin{equation}
	\widehat{v}^2 = - (v_1^2 + \cdots + v_n^2).
\end{equation}
More generally, if $v = (v_0, v_1, \ldots, v_n) \in {\bf R}^{n + 1}$,
\begin{equation}
	\widetilde{v} = v_0 + v_1 \, e_1 + \cdots + v_n \, e_n,
\end{equation}
and
\begin{equation}
	\widetilde{v}^* = v_0 - v_1 \, e_1 - \cdots - v_n \, e_n,
\end{equation}
then
\begin{equation}
	\widetilde{v} \, \widetilde{v}^* = \widetilde{v}^* \, \widetilde{v}
		= v_0^2 + v_1^2 + \cdots + v_n^2.
\end{equation}
Similarly, if $x = x_1 + x_2 \, i + x_3 \, j + x_4 \, k \in {\bf H}$,
where $x_1, x_2, x_3, x_4$ are real numbers, and we put
\begin{equation}
	x^* = x_1 - x_2 \, i - x_3 \, j - x_4 \, k,
\end{equation}
then
\begin{equation}
	x \, x^* = x^* \, x = x_1^2 + x_2^2 + x_3^2 + x_4^2.
\end{equation}

	If $f$ is a continuously-differentiable function on an open
set $U \subseteq {\bf R}^n$ with values in the Clifford algebra
$\mathcal{C}(n)$, then we say that $f$ is left or right Clifford
holomorphic if
\begin{equation}
	\sum_{l = 1}^n e_l \, \frac{\partial f}{\partial x_l} = 0
		\quad\hbox{or}\quad
	\sum_{l = 1}^n \frac{\partial f}{\partial x_l} \, e_l = 0
\end{equation}
on $U$, respectively.  Alternatively, let $f$ be a
continuously-differentiable function on an open set in ${\bf
R}^{n+1}$, where $x \in {\bf R}^{n + 1}$ has components $x_0, x_1,
\ldots, x_n$.  In this case, we get slightly different versions of
Clifford holomorphicity with the equations
\begin{equation}
	\sum_{l = 0}^n e_l \, \frac{\partial f}{\partial x_l} = 0, \quad
		\sum_{l = 0}^n \frac{\partial f}{\partial x_l} \, e_l = 0,
\end{equation}
where $e_0 = 1$.  There are also variants of these for the quaternions
using $i$, $j$, and $k$.  These are all \emph{Generalized
Cauchy--Riemann Systems} as in \cite{s-w}.

	Suppose that $f$ is a continuously-differentiable function on
an open set $U$ in ${\bf R}^n$ with values in ${\bf R}^p$ for some $p
\ge 1$, and that $x$ is an element of $U$ and $v$ is a unit vector in
${\bf R}^n$.  If $f$ satisfies an equation at $x$ like those described
in the previous paragraph, then the directional derivative of $f$ at
$x$ in the direction of $v$ can be expressed as a linear combination
of the directional derivatives of $f$ at $x$ in the directions
orthogonal to $v$.  For example, if $f$ is a left or right Clifford
holomorphic function, then one can check this by multiplying the
corresponding differential equation on the left or right by
$\widehat{v}$, respectively.  Let us say that $f$ is $k$-restricted
for some $k \ge 1$ if for every $x \in U$ and every hyperplane $H
\subseteq {\bf R}^n$, the norm of the differential of $f$ at $x$ is
less than or equal to $k$ times the norm of the restriction of the
differential of $f$ at $x$ to $H$.  In each of the cases discussed in
the previous paragraph, it follows that $f$ is $k$-restricted for a
fixed $k$.

	The differential of a real-valued function automatically
vanishes on a hyperplane at each point.  Hence a real-valued
$k$-restricted function on a connected open set is constant.  When $p
= n = 2$, the property of being $k$-restricted is very close to
quasiregularity.  A key difference is that quasiregularity includes a
condition of nonnegative orientation.

	Let us say that a linear mapping $A : {\bf R}^n \to {\bf R}^p$
is $k$-restricted if the norm of $A$ is less than or equal to $k$
times the norm of the restriction of $A$ to any hyperplane in ${\bf
R}^n$.  Equivalently, $A$ is $k$-restricted if the norm of $A$ is less
than or equal to $k$ times the norm of $A + B$ for every linear
mapping $B : {\bf R}^n \to {\bf R}^p$ with rank one.  This is also the
same as saying that the first singular value of $A$ is less than or
equal to $k$ times the second singular value.  Thus a
continuously-differentiable mapping is $k$-restricted if and only if
its differential is $k$-restricted at each point, which is exactly the
condition required for the arguments in \cite{s-w} for improved
subharmonicity properties of norms of vector-valued harmonic
functions.  It follows that $f : U \to {\bf R}^p$ is $k$-restricted if
and only if the norm of the differential of $f$ at any $x \in U$ is
less than or equal to $k$ times the norm of the differential of $f + a
\, \phi$ for every $a \in {\bf R}^p$ and continuously-differentiable
real-valued function $\phi$ on $U$.

	As an extension of quaternionic and Clifford analysis, one
could replace the usual partial derivatives in the coordinate
directions with vector fields with smooth coefficients.  The number of
vector fields could even be less than the dimension of the space, in
which event one might ask that the vector fields satisfy the
H\"ormander condition that they and their commutators span the tangent
space at each point.  This would imply that functions with vanishing
derivatives in the directions of the vector fields are locally
constant in particular.  Note that solutions of tangential
Cauchy--Riemann equations correspond to special classes of
quaternionic and Clifford holomorphic functions associated to suitable
vector fields, at least locally, just as for holomorphic functions and
quaternionic and Clifford analysis in the classical case.

	One can also consider versions of quaternionic and Clifford
analysis on metric spaces.  Since products of quaternionic or Clifford
holomorphic functions are not normally holomorphic even on Euclidean
spaces, abstract approaches based on algebras of functions do not work
as in the complex case.  One might look at $k$-restrictedness of a
function $f$ on a metric space in terms of comparing local Lipschitz
or Sobolev constants for $f$ with their counterparts for $f + a \,
\phi$ when $a$ is a constant vector and and $\phi$ is real-valued.

\section{Spaces with Poincar\'e inequalities}
\setcounter{equation}{0}

	If $B$ is a ball of radius $R > 0$ in ${\bf R}^n$, $1 \le p
< \infty$, and $f$ is a real-valued function on $B$, then
\begin{equation}
	\Big(\frac{1}{|B|} \, \int_B |f(x) - f_B| \, dx \Big)^{1/p}
\le C(n) \, R \, \Big(\frac{1}{|B|} \, \int_B |\nabla f(x)|^p \, dx\Big)^{1/p},
\end{equation}
where $|B|$ denotes the volume of $B$ and $f_B$ is the average of $f$
on $B$.  One might as well suppose that $f$ is continuously
differentiable on $B$, although the inequality also works when $f$ is
a locally integrable function on $B$ with distributional first
derivatives in $L^p(B)$.  The limiting case $p = \infty$ corresponds
to the statement that a Lipschitz condition is implied by a bound for
the gradient.

	Juha Heinonen and I posed some questions in \cite{h-s} about
whether suitable versions of these classical Poincar\'e inequalities
on other spaces would imply that the spaces enjoy some sort of
approximately Euclidean or sub-Riemannian structure.  These questions
were answered negatively by remarkable examples of Bourdon and Pajot
\cite{b-p} and Laakso \cite{l1}.  Perhaps it is better to say that
they answered positively the question of whether there could be a lot
of spaces of this type.  In particular, there are spaces of this type
with any Hausdorff dimension greater than or equal to $1$, and every
such space with at least two elements has Hausdorff dimension greater
than or equal to $1$ because of connectedness.

	I would like to suggest that there are positive results along
the lines of the previous questions with additional hypotheses.  There
is a nice theorem of Berestovskii and Vershik \cite{b-v} concerning
sub-Riemannian geometry of metric spaces under somewhat different
conditions, and one may be able to build on their approach.  Cheeger's
work \cite{che} on differentiability of Lipschitz functions almost
everywhere on spaces with Poincar\'e inequalities ought to be an
important step in this direction as well.  It may be relatively easy
to deal with spaces on which there is sufficient ``calculus'', and
there can be different amounts of structure corresponding to different
degrees of calculus.

	It can be helpful to look at nilpotent Lie groups and
sub-Riemannian spaces more closely in order to understand the general
situation better.  One can also simply start with a connected smooth
manifold $M$ and some smooth vector fields $X_1, \ldots, X_n$ on $M$
which satisfy the H\"ormander condition that the tangent space of $M$
at every point is spanned by the $X_\ell$'s and their successive Lie
brackets.  The smooth functions on $M$ as a smooth manifold would be
the same as the smooth functions on $M$ with respect to $X_1, \ldots,
X_n$, but the two structures can measure smoothness in very different
ways.  A vector field on $M$ is a first-order differential operator in
the usual sense but may be considered as an operator of higher order
with respect to the $X_\ell$'s.  A vector field which can be expressed
as the bracket of $r$ of the $X_\ell$'s in some way would typically be
considered a differential operator of order $r$ with respect to the
$X_\ell$'s.  In particular, one might start with a smooth distribution
$\mathcal{L}$ of linear subspaces of the tangent spaces of $M$, which
contain the $X_\ell$'s and are spanned by them at every point.  The
H\"ormander condition is then a maximal non-integrability condition
for $\mathcal{L}$, at least if $\mathcal{L}$ consists of proper
subspaces of the tangent spaces of $M$.

	However, that brackets of the $X_\ell$'s can be defined at all
can be considered as an important integrability condition for the
corresponding sub-Riemannian space.  To have any nontrivial vector
fields on a metric space at all is already quite significant, in the
sense of first-order differential operators acting on Lipschitz
functions as in \cite{nw3}, for instance.  Even if there are a lot of
vector fields on metric spaces with Poincar\'e inequalities by
\cite{che}, it may not be clear how to deal with their brackets.  In
the context of complex-analytic metric spaces, it would be interesting
to know whether brackets of complex vector fields of $\partial /
\partial \overline{z}$ type are of the same type.  This is the
classical integrability condition for an almost-complex structure on a
smooth manifold.

	There are classical results about integrating vector fields to
get nice mappings on manifolds.  Extra compatibility conditions are
required on sub-Riemannian spaces to get mappings which respect the
geometry in appropriate ways.  Even for nilpotent Lie groups with
sub-Riemannian structures that are invariant under left or right
translations by definition, there may only be finite-dimensional
families of mappings with suitable regularity.  In some cases there
may be no reason for a metric space with Poincar\'e inequalities to
have any nontrivial continuous families of mappings which respect the
geometry or the topology.  On complex-analytic metric spaces, it would
be interesting to consider holomorphic vector fields and the
possibility of integrating them to get holomorphic mappings.

	As another basis for comparison, suppose that $M$ is a smooth
manifold and that $V$ is a continuous vector field on $M$ which may
not be smooth.  There are still results about existence of integral
curves for $V$ in $M$, but uniqueness might not hold without more
information about the regularity of $V$.  Uniqueness can also fail for
smooth vector fields on singular spaces.  Similarly, let $L$ be a
Laakso space, with a projection from the Cartesian product of a Cantor
set $C$ and the unit interval $I$ onto $L$.  One can follow the
standard vector field on $I$ and move in the positive direction at
unit speed, and do the same on each parallel copy of $I$ in $C \times
I$.  Because of the identifications between the copies of $I$ in $L$,
one loses uniqueness of the trajectories in $L$.  It seems interesting
to consider metric spaces with vector fields more broadly, including
constructions like Laakso's with different patterns of
identifications.  One might wish to use probability theory to treat
this type of branching, i.e., to follow a vector field with a
stochastic process.

	On Laakso's and related spaces, there are nice classes of
regular functions which are locally equivalent to smooth functions on
the unit interval and constant in the direction of the Cantor set.  A
regular vector field can be defined as a regular function times
ordinary differentiation in the direction of the unit interval.  A
regular vector field applied to a regular function is a regular
function, and the bracket of two regular vector fields makes sense and
is a regular vector field.  Even for regular functions, many of the
usual problems are still present.  The branching can take place on
larger regions.

	Suppose now that $M$ is a smooth manifold equipped with some
sort of sub-Riemannian structure.  If $V$ is a vector field on $M$
which is smooth with respect to the ordinary smooth structure on $M$,
then one can integrate $V$ to get smooth mappings on $M$ which are at
least continuous with respect to the sub-Riemannian geometry.  If $V$
is admissible for the sub-Riemannian structure, then the integral
curves for $V$ are automatically admissible.  If $V$ is admissible and
$[V, X]$ is admissible when $X$ is, then the mappings on $M$
associated to $V$ are compatible with the sub-Riemannian structure.
This is basically a regularity condition for $V$ relative to the
sub-Riemannian structure, analogous to the classical Lipschitz
condition for the coefficients of a vector field.

	Sometimes a vector field is obtained from the gradient of a
function.  This could be derived from a pairing between functions
which includes a pointwise pairing between their gradients, at least
implicitly.  Such a pairing might be positive and symmetric, like a
Riemannian metric, or antisymmetric, as for a symplectic structure.

	On a Cantor set represented as the Cartesian product of a
sequence of finite sets, there are a lot of transformations obtained
from permutations in the individual coordinates, including small
displacements from permutations in coordinates with large indices.  It
is not so easy to have nontrivial small displacements on some locally
connected spaces, because of intricate topological structure.  On
spaces with Poincar\'e inequalities, one can also try to study small
displacements in terms of vector fields, perhaps on associated tangent
objects.  At any rate, it seems interesting to look at group actions
on spaces with Poincar\'e or related inequalities.

	I would like to think of a metric space with Poincar\'e,
Sobolev, or similar inequalities as a kind of generalized manifold.
One can look at this as a variant of \emph{Cantor manifolds}
\cite{h-w}, which are spaces that are not disconnected by subsets of
topological codimension $\ge 2$.  This is especially close to
isoperimetric inequalities and other estimates of the measure of a set
in terms of the size of its boundary.

	A lot of analysis related to singular integral operators and
classical spaces of functions is available in the vast setting of
spaces of homogeneous type \cite{c-w-1, c-w-2}.  This includes Cantor
sets and snowflake spaces, which are important examples with many
applications, and which also have nonconstant functions with vanishing
gradient.  As a next step, one can try to integrate local Lipschitz
conditions to estimate the behavior of a function.  With Poincar\'e or
Sobolev inequalities, one has stronger forms of calculus involving
integrals of derivatives.

	Notions of generalized manifolds have been studied extensively
in algebraic topology.  After all, homology and cohomology are also
ways of doing ``calculus'' on broad classes of spaces.  One of the
simplest of these notions is that of polyhedral pseudomanifolds.  More
sophisticated theories deal with intermediate dimensions in the space.

	Even for topological manifolds, there can be different types
of structures with different versions of calculus.  A manifold may be
equipped with a smooth structure, for instance, or a piecewise-linear,
Lipschitz, or quasiconformal structure more generally.  It may be
represented as a polyhedron, which might or might not have
piecewise-linear local coordinates, or coordinates of moderate
complexity.  Using sub-Riemannian geometry, a manifold can be a
fractal and still have Poincar\'e and Sobolev inequalities.

	The apparent irregularities of some spaces with Poincar\'e or
Sobolev inequalities could be attributed to using classical geometry
instead of something like noncommutative geometry.  With
noncommutative geometry, one can try to avoid complicated patchworks
of interconnections and gain local or infinitesimal symmetries.
Spaces with some version of calculus are basically extensions of
smooth manifolds whether they are based on classical or noncommutative
geometry, and one might as well try to work with both.

	In Connes' theory \cite{cns}, commutators between singular
integral operators and multiplication operators are like derivatives
and the Dixmier trace corresponds to integration.  The Dixmier trace
is an asymptotic version of the trace that applies to slightly more
than the usual trace class operators and vanishes on trace class
operators.  These asymptotic properties are important for localization
in the theory.  One-dimensional spaces are somewhat exceptional
because of unusual regularity of commutator operators.

	Of course, one-dimensional sets are exceptional in more
classical ways as well.  Connectedness plays a large role in this.  A
nice historical survey related to curvature can be found in the
introduction to \cite{ppl}.  Some amazing discoveries about singular
integral operators and complex analysis are explained in \cite{m-m-v}.

	At any rate, the general area of analysis on metric spaces has
seen a lot of activity, and it seems to me that there is plenty of
room for more.

\end{document}